\pdfoutput=1
\documentclass[12pt]{amsart}

\usepackage{amssymb}

\usepackage{amsthm}

\usepackage{latexsym}
\usepackage{amsmath}
\usepackage{amssymb}
\usepackage{amsthm}
\usepackage[all]{xy}
\usepackage{hyperref}
\hypersetup{
    colorlinks=true,       
    linkcolor=blue,          
    citecolor=blue,        
    filecolor=blue,      
    urlcolor=blue           
}

\newtheorem{thm}{Theorem}
\newtheorem{prop}[thm]{Proposition}
\newtheorem{cor}[thm]{Corollary}

\newtheorem{rem}[thm]{Remark}

\begin{document}

\title{Hypersurfaces achieving the Homma-Kim bound}

\author{Andrea Luigi Tironi}

\date{\today}

\address{
Departamento de Matem\'atica,
Universidad de Concepci\'on,
Casilla 160-C,
Concepci\'on, Chile}
\email{atironi@udec.cl}


\thanks{During the preparation of this paper, the author was partially supported 
by Proyecto VRID N. 214.013.039-1.OIN and the Project Anillo ACT 1415 PIA CONICYT}

\maketitle

\begin{abstract}
Let $X^n$ be a hypersurface in $\mathbb{P}^{n+1}$ with $n\geq 2$
defined over a finite field. The main result of this note is the
classification, up to projective equivalence, of hypersurfaces $X^n$
as above without a linear component when the number of their
rational points achieves the Homma-Kim bound.

\end{abstract}

\section{Introduction}

In a series of papers \cite{HK1, HK2, HK3}, Homma and Kim settled
the Sziklai conjecture \cite{S} for plane curves. In particular,
as a consequence of their results one can deduce that for any
plane curve $C$ of degree $d$ over a finite field $\mathbb{F}_q$
of $q$ elements without $\mathbb{F}_q$-linear components, the
number $N_q(C)$ of $\mathbb{F}_q$-points of $C$ is bounded by
$N_q(C)\leq (d-1)q+2$ and equality holds if and only if $d=q=4$
and $C$ is projectively equivalent to the plane curve
{\small
$$
X_0^4+X_1^4+X_2^4+X_0^2X_1^2+X_1^2X_2^2+X_0^2X_2^2+X_0^2X_1X_2+
X_0X_1^2X_2+X_0X_1X_2^2=0.$$
}
\noindent In \cite{HK4} the authors establish also an upper bound
for the number $N_q(X^n)$ of $\mathbb{F}_q$-points of a hypersurface
$X^n\subset\mathbb{P}^{n+1}$ of degree $d$ and dimension $n\geq 2$
which is an analogous to their bound for a plane curve. Moreover,
they show that their bound
$$\Theta_n^{d,q}:=(d-1)q^n+dq^{n-1}+q^{n-2}+\dots+q+1$$
is the best one for hypersurfaces without $\mathbb{F}_q$-linear components that
is linear on their degrees, because, for each finite field, they
give three nonsingular surfaces of different degrees that reach
their bound.

In line with the above results for plane curves we characterize, up to projective equivalence, all the
hypersurfaces $X^n\subset\mathbb{P}^{n+1}$ without
$\mathbb{F}_q$-linear components which reach the Homma-Kim bound
$\Theta_n^{d,q}$ by proving the following

\medskip

\begin{thm}\label{main theorem}
Let $X^n\subset\mathbb{P}^{n+1}$ be a hypersurface of degree $d\geq
2$ and dimension $n\geq 2$ defined over $\mathbb{F}_q$ and without
$\mathbb{F}_q$-linear components. Then $N_q(X^n)\leq\Theta_n^{d,q}$
and equality holds if and only if $d\leq q+1$ and one of the
following possibilities occurs:
\begin{enumerate}
\item[$(1)$] $d=q+1$ and $X^n$ is a space filling hypersurface
$$(X_0,...,X_{n+1})\ A\ {}^t \! (X_0^q,...,X_{n+1}^q)=0,$$ where $A=\left( a_{ij}\right)_{i,j=1,...,n+2}$ is an
$(n+2)\times (n+2)$ matrix such that ${}^t \! A=-A$ and $a_{kk}=0$ for every $k=1,...,n+2$; moreover, $X^n$ is nonsingular if and only if $\det A\neq 0$;
\item[$(2)$]
$d=\sqrt{q}+1$ and $X^n$ is projectively equivalent to a cone over
the nonsingular Hermitian surface 
$$X_0^{\sqrt{q}+1}+X_1^{\sqrt{q}+1}+X_2^{\sqrt{q}+1}+X_3^{\sqrt{q}+1}=0;$$
\item[$(3)$] $d=2$ and $X^n$ is projectively equivalent to a quadric
hypersurface

\smallskip

$X_0(a_0X_0+...+a_{n+1}X_{n+1})+X_1(b_0X_0+...+b_{n+1}X_{n+1})=0,$

\smallskip

\noindent with $a_i,b_i\in\mathbb{F}_q$, which is nonsingular if and only if
$n=2$ and $$\det
\left(
\begin{array}{cc}
  a_2 & a_{3} \\
  b_2 & b_{3} \\
\end{array}
\right)\neq 0.$$
\end{enumerate}
\end{thm}

Finally, as a consequence of this result, in the
nonsingular case we deduce the following

\begin{cor}\label{Cor1}
Let $V\subset\mathbb{P}^{n+1}$ be a nonsingular hypersurface of
degree $d\geq 2$ and dimension $n\geq 2$ defined over
$\mathbb{F}_q$. Then $N_q(V)\leq \Theta_n^{d,q}$ and equality holds if and only if one of the following
cases occurs:
\begin{enumerate}
\item[$(a')$] $n$ is even and $V$ is a space filling hypersurface
$$(X_0,...,X_{n+1})\ A\ {}^t \! (X_0^q,...,X_{n+1}^q)=0,$$ where $A$ is an
$(n+2)\times (n+2)$ matrix such that ${}^t \! A=-A$, $\det A\neq 0$ and $a_{kk}=0$ for every $k=1,...,n+2$;

\smallskip

\item[$(b')$] $n=2$ and $V$ is projective equivalent to the Hermitian surface
$$X_0^{\sqrt{q}+1}+X_1^{\sqrt{q}+1}+X_2^{\sqrt{q}+1}+X_3^{\sqrt{q}+1}=0;$$

\item[$(c')$] $n=2$ and $V$ is projective equivalent to the hyperbolic quadric surface
$$X_0X_1-X_2X_3=0.$$
\end{enumerate}
\end{cor}

\medskip

\noindent\textit{Note}. After this paper was written, the author
discovered that in \cite{HK6} Homma and Kim characterize all the
surfaces in $\mathbb{P}^3$ without $\mathbb{F}_q$-lineal
components which reach their bound $\Theta_2^{d,q}$ by proving
Theorem \ref{main theorem} for $n=2$. Although the approach in
\cite{HK6} seems only slightly different to that of this note, for
completeness we preferred do not omit the proof of this case.

\medskip

\noindent\textit{Acknowledgements}. The author would like to
thank Antonio Laface for his kind and constant
encouragement and for many interesting discussions about algebraic
geometry and the referee for useful remarks which allow him to 
improve the presentation of this paper.

\section{Notation and preliminary results}

Let $X^n\subset\mathbb{P}^{n+1}$ be a hypersurface of degree $d\geq
2$ and dimension $n\geq 2$ defined over a finite field
$\mathbb{F}_q$ of $q$ elements, with $q=p^s$ for some prime number
$p$ and an integer $s\in\mathbb{Z}_{\geq 1}$. Denote by $N_q(W)$
the number of $\mathbb{F}_q$-points of a set
$W\subseteq\mathbb{P}^{n+1}$, by $W(\mathbb{F}_q)$ the set of
$\mathbb{F}_q$-points of $W$ and by $W^{\star}$ the set of all
hyperplanes $\mathbb{P}^n\subset\mathbb{P}^{n+1}$ containing the
linear span $\langle W \rangle$ of $W$. Recall that for any
$N\in\mathbb{Z}_{\geq 1}$ we have
$$N_q(\mathbb{P}^{N})=q^N+q^{N-1}+\dots +q+1.$$

In this section, we give some remarks and preliminary results
which will be useful to prove Theorem \ref{main theorem}.

\begin{rem}\label{Remark1}
By definition, $\Theta_n^{d,q}\leq N_q(\mathbb{P}^{n+1})$ if and
only if $d\leq q+1$.
\end{rem}

\begin{prop}\label{Prop1}
Let $X^n\subset\mathbb{P}^{n+1}$ be a hypersurface of degree $d\geq
2$ and dimension $n\geq 2$ defined over $\mathbb{F}_q$. If
$N_q(X^n)=\Theta_n^{d,q}$, then $d\leq q+1$ and equality holds if
and only if $N_q(X^n)=N_q(\mathbb{P}^{n+1})$, i.e. $X^n$ is a space
filling hypersurface of $\mathbb{P}^{n+1}$.
\end{prop}

\noindent\textit{Proof.} Since $X^n\subset\mathbb{P}^{n+1}$ and $N_q(X^n)=\Theta_n^{d,q}$, we see that
$$\Theta_n^{d,q}=N_q(X^n)\leq N_q(\mathbb{P}^{n+1}).$$
Thus by Remark \ref{Remark1} we obtain that $d\leq q+1$. Moreover,
we get $d=q+1$ if and only if $N_q(X^n)=N_q(\mathbb{P}^{n+1})$,
which is equivalent to say that $X^n$ is a space filling
hypersurface of $\mathbb{P}^{n+1}$. \hfill{$\square$}

\begin{rem}\label{Remark2}
From the proof of Case $(2)$ of {\em Theorem 1.2} in {\em
\cite[$\S 3.1$]{HK4}} it follows that if $d\leq q$ and $X^n$ does
not contain a linear $\mathbb{P}^{n-1}$ then
$N_q(X^n)<\Theta_n^{d,q}$. Therefore, if $N_q(X^n)=\Theta_n^{d,q}$ and
$d\leq q$ then there exists at least a linear space
$\mathbb{P}^{n-1}$ in $X^n$.
\end{rem}

\begin{thm}[\cite{Se}]\label{S-S-S}
Let $X^n\subset\mathbb{P}^{n+1}$ be a hypersurface of degree $d$
defined over $\mathbb{F}_q$. Then $N_q(X^n)\leq
dq^{n}+q^{n-1}+...+q+1$. Furthermore, if $d\leq q$ then equality
holds if and only if $X^n$ is a union of $d$ hyperplanes over
$\mathbb{F}_q$ that contain a common linear subspace of
codimension $2$ in $\mathbb{P}^{n+1}$.
\end{thm}

\begin{prop}\label{Prop2}
Let $X^n\subset\mathbb{P}^{n+1}$ be a hypersurface 
of degree $d\geq 2$ and dimension $n\geq 2$ defined over
$\mathbb{F}_q$. If $N_q(X^n)=\Theta_n^{d,q}$ and $d\leq q$, then $X^n$
is covered by linear spaces of dimension $n-1$.
\end{prop}

\noindent\textit{Proof.} By Remark \ref{Remark2} we know that
there exists a linear $L:=\mathbb{P}^{n-1}$ in $X^n$. Since
$N_q(L\cap X^n)=N_q(L):=d_2$, we have
$$N_q(X^n)=\left[ \sum_{H_i\in L^{\star}}\left( N_q(H_i\cap X^n)-d_2\right)\right] +d_2=$$
$$=\left[ \sum_{H_i\in L^{\star}} N_q(H_i\cap X^n)\right]-d_2q=$$
$$=\left[ \sum_{H_i\in L^{\star}} N_q(H_i\cap X^n)\right]-(q^{n-1}+q^{n-2}+\dots +q+1)q,$$
where $L^{\star}$ denotes the set of hyperplanes of
$\mathbb{P}^{n+1}$ containing $L$. Then we get
$$\left[ \sum_{H_i\in L^{\star}} N_q(H_i\cap X^n)\right]=(q+1)(dq^{n-1}+q^{n-2}+ \dots +q+1)$$
and by Theorem \ref{S-S-S} we deduce that $X^n\cap
H_i=\cup_{k=1}^{d}\mathbb{P}_{i,k}^{n-1}$ for every $H_i\in
L^{\star}$. Thus we conclude that
$$X^n=X^n\cap\mathbb{P}^{n+1}=X^n\cap\left(\cup_{H_i\in L^{\star}}\ H_i\right)=\cup_{H_i\in L^{\star}}\ \left(X^n\cap H_i\right)=\cup_{i,k}\ \mathbb{P}_{i,k}^{n-1},$$ i.e. $X^n$ is covered by
linear $\mathbb{P}^{n-1}$'s. \hfill{$\square$}

\begin{rem}\label{Remark3}
From the proof of {\em Proposition \ref{Prop2}} it follows that if
$X^n$ contain a linear $L=\mathbb{P}^{n-1}$ then for any $H\in
L^{\star}$ we have $X^n\cap H=\cup_{i=1}^{d} L_i$, where the
$L_i$'s are linear $\mathbb{P}^{n-1}$ such that, after renaming,
$L_1=L$ and all the $L_i$'s meet in a common subspace
$\mathbb{P}^{n-2}$.
\end{rem}

\begin{prop}\label{Prop3}
Let $X^n\subset\mathbb{P}^{n+1}$ be a hypersurface of degree $d\geq
2$ and dimension $n\geq 3$ defined over $\mathbb{F}_q$. 
If $N_q(X^n)=\Theta_n^{d,q}$ and $d\leq q$ then $X^n$ is singular,
that is, $\mathrm{Sing}(X^n)\neq\emptyset$.
\end{prop}

\noindent\textit{Proof.} By Remark \ref{Remark2} we know that
there exists a linear $L:=\mathbb{P}^{n-1}$ in $X^n$. Furthermore,
from Remark \ref{Remark3} we deduce that $X^n\cap
H_i=\cup_{k=1}^{d}\mathbb{P}_{i,k}^{n-1}$ for every $H_i\in
L^{\star}$. Thus consider two distinct hyperplanes $H_1,H_2\in
L^{\star}$ such that
$$X^n\cap H_1=\cup_{k=2}^{d}\mathbb{P}_{1,k}^{n-1}\cup L, \ \ X^n\cap H_2=\cup_{h=2}^{d}\mathbb{P}_{2,h}^{n-1}\cup L.$$
Since $d\geq 2$, up to renaming, we see that there are two
distinct linear spaces $\mathbb{P}_{1,2}^{n-1}\subset X^n\cap H_1$
and $\mathbb{P}_{2,2}^{n-1}\subset X^n\cap H_2$ such that
$\Lambda_i:=L\cap\mathbb{P}_{i,2}^{n-1}$ is a linear space
$\mathbb{P}^{n-2}$ in $L$ for $i=1,2$. Therefore we have
$$n-1=\dim L\geq\dim (\Lambda_1+\Lambda_2)=\dim \Lambda_1+\dim \Lambda_2-\dim (\Lambda_1\cap \Lambda_2),$$
i.e. $\dim (\Lambda_1\cap \Lambda_2)\geq n-3\geq 0$. This shows the existence
of at least a point $p\in \Lambda_1\cap \Lambda_2\subset H_1\cap H_2$ such
that $p\in L\subset X^n$ and the tangent space $T_pX^n$ contains $H_1$
and $H_2$, i.e. $H_1\cup H_2\subset T_pX^n$. This implies
$T_pX^n=\mathbb{P}^{n+1}$, i.e. $p\in X^n$ is a singular point of $X^n$.
\hfill{$\square$}

\begin{rem}\label{Remark4}
Let $Y^n\subset\mathbb{P}^{n+1}$ be a hypersurface of dimension
$n$ over $\mathbb{F}_q$. Suppose that $Y^n$ is a cone
$\mathcal{C}_{\mathbb{P}^k}(Y^{n-k-1})$ over a hypersurface
$Y^{n-k-1}\subset\mathbb{P}^{n-k}$ of dimension $n-k-1\geq 2$ with
vertex $\mathbb{P}^k$ such that $0\leq k\leq n-3$ and
$\mathbb{P}^k\cap\mathbb{P}^{n-k}=\emptyset$ in
$\mathbb{P}^{n+1}$. Then we have the following properties $($see
also {\em \cite[(3.1),(3.2)]{HK4}}$):$
\begin{enumerate}
\item[(a)] $d:=\deg Y^{n}=\deg Y^{n-k-1}$;

\item[(b)] $N_q(Y^{n})=q^{k+1}N_q(Y^{n-k-1})+N_q(\mathbb{P}^k)$;

\item[(c)] $N_q(Y^{n})=\Theta_n^{d,q}$ if and only if
$N_q(Y^{n-k-1})=\Theta_{n-k-1}^{d,q}$.

\end{enumerate}
\end{rem}

\begin{prop}\label{Prop4}
Let $X^n\subset\mathbb{P}^{n+1}$ be a hypersurface of degree $d\geq
2$ and dimension $n\geq 2$ defined over $\mathbb{F}_q$. If
$N_q(X^n)=\Theta_n^{d,q}$ and $d\leq q$, then $X^n$ is a reduced hypersurface and the set 
of singular points $\mathrm{Sing}(X^n)$ of $X^n$ is a proper closed subset of $X^n$,
that is, $\overline{\mathrm{Sing}(X^n)}=\mathrm{Sing}(X^n)\subsetneq X^n$.
\end{prop}

\noindent\textit{Proof.} Suppose that $X^n$ contains a nonreduced component.
Then $d':=\deg X_{red}^n < \deg X^n=d$ and $N_q(X^n)=N_q(X_{red}^n)$, where $X_{red}^n$ is the reduced
hypersurface obtained from $X^n$. By Theorem \ref{S-S-S} we have
$$\Theta_n^{d,q}=N_q(X^n)=N_q(X_{red}^n)\leq d'q^n+q^{n-1}+...+1\leq (d-1)q^n+q^{n-1}+...+1,$$ but this 
gives $d\leq 1$, a contradiction. Thus $X^n$ is a reduced hypersurface and the statement follows from
\cite[Theorem 5.3]{H}. \hfill{$\square$}

\section{Proof of Theorem $\ref{main theorem}$}

Let us recall that from Proposition \ref{Prop1} it follows that
$d\leq q+1$, with equality if and only if $X^n$ is a space filling
hypersurface with $N_q(X^n)=N_q(\mathbb{P}^{n+1})$. Moreover, if
$d=2$ then $X^n$ is a quadric hypersurface of $\mathbb{P}^{n+1}$.

\noindent Therefore, without loss of generality, we can assume that $$3\leq
d \leq q.$$

\noindent Furthermore, by Proposition \ref{Prop4}, let $p\in X^n$ be a
nonsingular point.

\bigskip

\noindent First of all, assume that $n=2$, i.e.
$X^2\subset\mathbb{P}^3$ is a surface. By \cite[(3.8)]{HK5} there
exists a plane $H=\mathbb{P}^2$ such that $H\cap
X^2=l_1\cup\dots\cup l_d$ with $\cap_{i=1}^{d}l_i=\{p\}$, where the
$l_j$'s are $d$ lines $\mathbb{P}^1$ in $X^2$. Take a line $L\subset
H$ such that $L\cap X^2=\{p\}$ and consider a pencil of planes
$H_i\in L^{\star}$ for $i=1,...,q+1$ with $H_{q+1}=H$. Put
$\Gamma_i:=H_i\cap X^2\subset H_i=\mathbb{P}^2$ for $i=1,...,q+1$
and note that $\Gamma_{q+1}=l_1\cup\dots\cup l_{d}$. Furthermore,
observe that
$$\Theta_2^{d,q}=N_q(X^2)=\left[\sum_{i=1}^q\left(N_q(\Gamma_i)-1\right)\right]+N_q(\Gamma_{q+1})=$$
$$=\left[\sum_{i=1}^qN_q(\Gamma_i)\right]-q+dq+1,$$ i.e.
\begin{equation}
\sum_{i=1}^qN_q(\Gamma_i)=\left[(d-1)q+1\right]q.\tag{*}
\end{equation}

\bigskip

\noindent\textit{Claim 1.} $N_q(\Gamma_i)=(d-1)q+1$ for every
$i=1,...,q$.

\smallskip

\noindent Note that $T_pX^2=H$. Moreover, since $p$ is a nonsingular
point of $X^2$, we see that all the $\Gamma_i$'s are plane curves
without linear components. Thus if either $q\neq 4$ or $d\neq 4$, then
$N_q(\Gamma_i)\leq (d-1)q+1$ for any $i=1,...,q$. On the other
hand, if $q=4=d$ then $N_q(\Gamma_i)\leq (d-1)q+2=14$. Suppose,
after renaming, that $N_q(\Gamma_1)= 14$ and take a bitangent line
$l'\subset H_1$ to $\Gamma_1$. Since $d\geq 3$, by
\cite[(3.6)]{HK5} we deduce that $\hat{H}_i\cap X^n$ has not linear
components for any $\hat{H}_i\in(l')^{\star}$. Thus we get
$$65=N_q(X^2)=N_q(\Gamma_1)+\sum_{H_1\neq\hat{H}_i\in(l')^{\star}}
\left[N_q(\hat{H}_i\cap X^2)-2\right]\leq $$
$$\leq 14+4(14-2)=62,$$ but this
is a contradiction. Thus $N_q(\Gamma_i)\leq (d-1)q+1$ for any
$i=1,...,q$ and the statement follows from (*). \hfill{Q.E.D.}

\bigskip

\noindent From \cite[$\S 2$]{HK2} we know that any $\Gamma_i$ is
absolutely irreducible and
$$N_q(\mathrm{Sing}(\Gamma_i)\cap\Gamma_i)=0\ \ \mathrm{ for }\ i=1,...,q.$$

\bigskip

\noindent\textit{Claim 2}. $N_q(l\cap\Gamma_i)\in\{0,1,d\}$ for
any line $l\subset H_i=\mathbb{P}^2$ with $i=1,...,q$.

\smallskip

\noindent Assume that there exists a line $l'\subset
H_i=\mathbb{P}^2$ such that $2\leq N_q(l'\cap\Gamma_i)\leq d-1$ for some $i\in \{1,...,q \}$.
Consider a pencil $L_j\in (l')^{\star}$ of planes
$L_j=\mathbb{P}^2$ which contain the line $l'$. Then by
\cite[(3.6)]{HK5} there is not a plane $L_j\in (l')^{\star}$ such
that $L_j\cap X^2$ contains a line. Thus by \cite{HK3} we conclude
that
\begin{equation*}
\Theta_2^{d,q}=N_q(X^2)=N_q(\Gamma_i)+\sum_{H_i\neq L_j\in
{l'}^{\star}}\left[ N_q(L_j\cap X^2)-N_q(l'\cap\Gamma_i)\right]\leq
\end{equation*}
\begin{equation*}
\leq
(d-1)q+1+\sum_{j=1}^{q}\left[(d-1)q+2-N_q(l'\cap\Gamma_i)\right]\leq
(d-1)q^2+dq,
\end{equation*}
but this is a contradiction. \hfill{Q.E.D.}

\bigskip

\noindent Therefore, by counting the number $\delta_i$ of distinct
lines $l\subset H_i=\mathbb{P}^2$ such that $N_q(l\cap
\Gamma_i)>0$ we have
$$\frac{\binom{N_q(\Gamma_i)}{2}}{\binom{d}{2}}+N_q(\Gamma_i)=\delta_i\leq N_q({\mathbb{P}^2}^*)=N_q(\mathbb{P}^2)=
q^2+q+1,$$ where ${\mathbb{P}^2}^*$ is the dual of $\mathbb{P}^2$.
This gives $\frac{(d-1)q+1}{d}+d-1\leq q+1,$ that is,
$d^2-2d-(q-1)\leq 0$. Hence $3\leq d\leq\sqrt{q}+1$.

If $d=\sqrt{q}+1$ then from \cite{HK5} we know that $X^2$ is a
nonsingular Hermitian surface. Finally, suppose that $3\leq
d<\sqrt{q}+1$. By \cite{AP} we see that
$$N_q(C)\leq q+1+2p_C\sqrt{q}<(d-1)q+1$$ for any absolutely irreducible plane curve $C$, where
$p_C=\frac{(d-1)(d-2)}{2}$ is the arithmetic genus of $C$.
Therefore, we obtain that
$$\left[(d-1)q+1\right]q=\sum_{i=1}^qN_q(\Gamma_i)\leq [(d-1)q]q,$$ which gives again a numerical contradiction.

\medskip

\noindent All the above arguments show that if $n=2$ then
$d\in\{2,\sqrt{q}+1,q+1\}$.

\medskip

Assume now that $n\geq 3$. Since $d\leq q$, from Proposition \ref{Prop2} we know that $X^n$ is covered by linear
$\mathbb{P}^{n-1}$'s and by Proposition \ref{Prop4} we can consider a linear space $L=\mathbb{P}^{n-1}$
which contains a nonsingular point $p\in X^n$.

Let $H_0$ be a linear $\mathbb{P}^n$ such that $L\subset H_0$.
By Remark \ref{Remark3} we get $X^n\cap H_0=\cup_{i=1}^{d}L_i$ with $L_i=\mathbb{P}^{n-1}$,
$L=L_1$ and $\cap_{i=1}^{d}L_i:=\Lambda=\mathbb{P}^{n-2}\subset
L$. Then we have two possibilities: (i) $p\in\Lambda$; \ (ii) $p\notin\Lambda$.

In case (i), since $n-2\geq 1$, we see that there exists a linear space $L'=\mathbb{P}^{n-1}\subseteq H_0$
such that $p\in L'$, $\dim (L'\cap L_i)=n-2$ for $i=1,...,d$ and $\dim (L'\cap \Lambda)=n-3\geq 0$.
Note that $L'$ is not contained in $X^n$. Consider now a hyperplane $H_i\neq H_0$ which contains $L'$. If
there exists a linear $\mathbb{P}_{i_1}^{n-1}\subseteq X^n\cap H_i$,
then there exist $d$ linear spaces
$\mathbb{P}_{i_k}^{n-1}$ such that $X^n\cap
H_i=\cup_{k=1}^{d}\mathbb{P}_{i_k}^{n-1}$. Since 
$$L'\cap X^n=L'\cap (H_0\cap X^n)=L'\cap (\cup_{i=1}^{d}L_i)=\cup_{i=1}^{d}(L'\cap L_i)$$
and $p\in L'\cap\Lambda=\cap_{i=1}^{d}(L'\cap L_i)$, we deduce that 
$p\in\cap_{k=1}^{d}\mathbb{P}_{i_k}^{n-1}$. Thus the tangent space $T_pX^n$ contains $L_1$
, $L_2$ and $\mathbb{P}_{i_j}^{n-1}$ for some $j=1,...,d$ such that $\mathbb{P}^{n+1}=\langle L_1\cup L_2\cup \mathbb{P}_{i_j}^{n-1} \rangle$. 
This implies that $\mathbb{P}^{n+1}\subseteq T_pX^n$, i.e. $p\in X^n$ is a singular point of $X^n$, but this is a contradiction.
So $X^n\cap H_i$ has no $\mathbb{F}_q$-linear components for any hyperplane $H_i\neq H_0$ which contains $L'$. 

Assume now we are in case (ii). Take $L'=\mathbb{P}^{n-1}\subseteq H_0$ such that $p\in L'$ and $L'$ is as in case (i).
Consider now a hyperplane $H_i\neq H_0$ which contains $L'$. 
If there exists a linear $\mathbb{P}_{i_1}^{n-1}\subseteq X^n\cap H_i$,
then there exist $d$ linear spaces 
$\mathbb{P}_{i_k}^{n-1}$ such that $X^n\cap
H_i=\cup_{k=1}^{d}\mathbb{P}_{i_k}^{n-1}$ and
$p\in\mathbb{P}_{i_j}^{n-1}$ for some $j=1,...,d$. Since $L_1\cap L' = \mathbb{P}^{n-2} = L' \cap \mathbb{P}_{i_j}^{n-1}$
and $L_1\neq \mathbb{P}_{i_j}^{n-1}$, we see that 
$$n-2\geq\dim (L_1\cap \mathbb{P}_{i_j}^{n-1}) \geq \dim (L_1\cap \mathbb{P}_{i_j}^{n-1}\cap L') =n-2,$$
i.e. $\dim (L_1\cap \mathbb{P}_{i_j}^{n-1})=n-2$. Therefore, by taking $H_0:=\langle L_1\cup \mathbb{P}_{i_j}^{n-1} \rangle$, we can repeat all the above arguments
and this allows us to get to case (i) again. Thus it follows that $X^n\cap H_i$ has no $\mathbb{F}_q$-linear components for any hyperplane $H_i\neq H_0$ which contains $L'$. 

In anyway, the above arguments show that in both cases
we can suppose that there exists a linear space $L'=\mathbb{P}^{n-1}\subset H_0$ not contained in $X^n$
such that $p\in L'$, $X^n\cap H_0$ is the union of $\mathbb{P}^{n-1}$'s meeting in a 
common linear subspace $\mathbb{P}^{n-2}$ and $X^n\cap H_i$ has no $\mathbb{F}_q$-linear components for any hyperplane $H_i\neq H_0$ which contains $L'$. 
This gives
$$\Theta_n^{d,q}=N_q(X^n)=\sum_{i=1}^{q}\left[N_q(X^n\cap
H_i)-N_q(\cup_{k=1}^{d}L'\cap L_k)\right]+N_q(X^n\cap H_0)=$$
$$=\left[\sum_{i=1}^{q}N_q(X^n\cap H_i)\right]-q(dq^{n-2}+q^{n-3}+...+1)+(dq^{n-1}+q^{n-2}+...+1),$$
i.e.
$$\sum_{i=1}^{q}\left[\Theta_{n-1}^{d,q}-N_q(X^n\cap
H_i)\right]=0.$$ Since $N_q(X^n\cap H_i)\leq\Theta_{n-1}^{d,q}$ for
every $i=1,...,q$, we get $N_q(X^n\cap H_i)=\Theta_{n-1}^{d,q}$. By
an inductive argument, we conclude that $d\in\{2,\sqrt{q}+1,q+1\}$
for $n\geq 2$.

\bigskip

Finally, the next three results 
allow us to conclude the proof of Theorem
\ref{main theorem}.

\begin{prop}[$d=2$]
Let $X^n\subset\mathbb{P}^{n+1}$ be a hypersurface of degree two
without $\mathbb{F}_q$-linear components and such that
$N_q(X^n)=\Theta_n^{2,q}$. If $n\geq 2$ then $X^n$ is projectively equivalent to a quadric
hypersurface
$$X_0(a_0X_0+...+a_{n+1}X_{n+1})+X_1(b_0X_0+...+b_{n+1}X_{n+1})=0,$$
with $a_i,b_i\in\mathbb{F}_q$, which is nonsingular if and only if
$n=2$ and $$\det
\left( \begin{array}{cc}
  a_2 & a_{3} \\
  b_2 & b_{3} \\
\end{array}
\right)\neq 0.$$
\end{prop}

\noindent\textit{Proof}. From Remark \ref{Remark2} it
follows that there exists a linear subspace $L:=\mathbb{P}^{n-1}\subset X^n$.
Thus, after a change of coordinates, we can suppose that $L=\{ X_0=X_1=0 \}$ and
$$X^n\ :\quad X_0L_0+X_1L_1=0,$$ where $L_i:=L_i(X_0,...,X_{n+1})$ are homogeneous polynomial of
degree one. Write $L_i:=\sum_{j=0}^{n+1}\alpha_j^iX_j$ for some
$\alpha_j^i\in\mathbb{F}_q$ and $F:=X_0L_0+X_1L_1$. Note that
$X^n=\{ F=0 \}$ and
$$\frac{\partial \vec{F}}{\partial X}:=
\left(
\begin{array}{c}
\frac{\partial F}{\partial X_0} \\[3pt]
\frac{\partial F}{\partial X_1} \\[3pt]
\frac{\partial F}{\partial X_2} \\[3pt]
\vdots \\[3pt]
\frac{\partial F}{\partial X_{n+1}}
\end{array}\right)=A\cdot 
\left(
\begin{array}{c}
X_0 \\[3pt]
X_1 \\[3pt]
X_2 \\[3pt]
\vdots \\[3pt]
X_{n+1}
\end{array}\right)=:A\cdot {}^t \! \vec{x},
$$
where $A$ is the following
$(n+2)\times (n+2)$ matrix
$$
\left(
\begin{array}{ccccc}
  2\alpha_0^0 & \alpha_1^0+\alpha_0^1 & \alpha_2^0 & \dots & \alpha_{n+1}^0 \\
  \alpha_0^1+\alpha_1^0 & 2\alpha_1^1 & \alpha_2^1 & \dots & \alpha_{n+1}^1 \\
  \alpha_2^0 & \alpha_2^1 & 0 & \dots & 0 \\
  \alpha_3^0 & \alpha_3^1 & 0 & \dots & 0 \\
  \vdots & \vdots & \vdots &  & \vdots \\
  \alpha_{n+1}^0 & \alpha_{n+1}^1 & 0 & \dots & 0 \\
\end{array}
\right).$$ 

\medskip

\noindent\textit{Claim 3}. \ $X^n$ is singular $\iff$ \ $\det A=0$.

\smallskip

\noindent Assume that $\det A=0$. Then there exists $\vec{y}\neq\vec{0}$ such that
$A\cdot {}^t \! \vec{y}=\vec{0}$. Write $\vec{y}:=(y_0,y_1,y_2,...,y_{n+1})$.
If $(y_0,y_1)\neq (0,0)$, then 
$0=\det\left(
\begin{array}{cc}
\alpha_2^0 & \alpha_2^1 \\
  \alpha_3^0 & \alpha_3^1 
  \end{array}
\right)=\det\left(
\begin{array}{cc}
\alpha_2^0 & \alpha_3^0 \\
  \alpha_2^1 & \alpha_3^1 
  \end{array}
\right)$. This shows that there exists $(z,w)\neq (0,0)$ such that
$\left(
\begin{array}{cc}
\alpha_2^0 & \alpha_3^0 \\
  \alpha_2^1 & \alpha_3^1 
  \end{array}
\right)\cdot \left(
\begin{array}{c}
z \\
w 
\end{array}
\right)=\left(
\begin{array}{c}
0 \\
0 
\end{array}
\right)$. Thus the vector $\vec{v}:=(0,0,z,w,0,...,0)$
is a nontrivial solution of $F=0$ and $\frac{\partial \vec{F}}{\partial X}(\vec{v})=A\cdot {}^t \! \vec{v}=\vec{0}$,
i.e. $X^n$ is singular. On the other hand, if $(y_0,y_1)=(0,0)$
then $\vec{y}\neq\vec{0}$ is a solution of $F=0$ and $\frac{\partial \vec{F}}{\partial X}(\vec{y})=A\cdot {}^t \! \vec{y}=\vec{0}$,
i.e. $X^n$ is singular again. 

Suppose now that $X^n$ is singular.
Then there exists $\vec{x}\neq\vec{0}$ such that
$A\cdot {}^t \! \vec{x}=\vec{0}$ and this implies that necessarily $\det A=0$.\hfill Q.E.D.

\bigskip

\noindent If $n\geq 3$, i.e. $n+1\geq 4$, then the last three rows of $A$ are linearly dependent, i.e. $\det
A=0$, and from \textit{Claim 3} we deduce that $X^n$ is singular. 
On the other hand, if $n=2$ then $X^2$ is given by
$$X_0(\alpha_0^0X_0+...+\alpha_{3}^0X_{3})+X_1(\alpha_0^1X_0+...+\alpha_{3}^1X_{3})=0.$$
Hence by \textit{Claim 3} we see that 
$$X^2\subset\mathbb{P}^3\ \mathrm{ is\ singular } \iff\det
\left(
\begin{array}{cc}
  \alpha_2^0 & \alpha_3^0 \\
  \alpha_2^1 & \alpha_3^1 \\
\end{array}
\right)= 0$$ and this proves the last part of the statement.
\hfill $\square$

\bigskip

\begin{prop}[$d=\sqrt{q}+1$]
Let $X^n\subset\mathbb{P}^{n+1}$ be a hypersurface of degree
$\sqrt{q}+1$ without $\mathbb{F}_q$-linear components and such that
$N_q(X^n)=\Theta_n^{\sqrt{q}+1,q}$. If $n\geq 2$ then $X^n$ is projectively equivalent to 
a cone over the nonsingular Hermitian surface
$$X_0^{\sqrt{q}+1}+X_1^{\sqrt{q}+1}+X_2^{\sqrt{q}+1}+X_3^{\sqrt{q}+1}=0.$$
\end{prop}

\noindent\textit{Proof}. If $n=2$ then from \cite{HK5} we know that $X^2$ is a nonsingular Hermitian surface $S$
in $\mathbb{P}^3$. Assume now that $n\geq 3$. 

By induction, suppose that any 
hypersurface $X^{n-1}\subset\mathbb{P}^{n}$ of degree $\sqrt{q}+1$ without $\mathbb{F}_q$-linear components 
and such that $N_q(X^{n-1})=\Theta_{n-1}^{\sqrt{q}+1,q}$ is 

\medskip

\noindent $(*)$ \ a cone 
over a nonsingular Hermitian surface $S$ and 
projectively equivalent to $\{ f=0\}\subset\mathbb{P}^{n}$, where
$$f:=(X_0,...,X_{n})\ A\ {}^t \! (X_0^{\sqrt{q}},...,X_{n}^{\sqrt{q}})$$
and $A$ is an $(n+1)\times (n+1)$ matrix such that ${}^t \! A=A^{(\sqrt{q})}$.

\medskip

From the proof of Theorem \ref{main theorem} in $\S 3$, we know
that there exist a linear $L_{11}=\mathbb{P}^{n-1}\subset X^n$, a
hyperplane $H_1=\mathbb{P}^n$ such that $L_{11}\subset H_1$, $X^n\cap
H_1=\cup_{j=1}^{\sqrt{q}+1}L_{1j}$ with $L_{1j}=\mathbb{P}^{n-1}$,
$\Lambda_1:=\cap_{j=1}^{\sqrt{q}+1}L_{1j}=\mathbb{P}^{n-2}$ and a
linear $L'=\mathbb{P}^{n-1}\subset H_1$ such that $L'\nsubseteq
X^n$, $X^n\cap H':=X^{n-1}$ has no linear
$\mathbb{F}_q$-components for any hyperplane $H'\neq H_1$ which contains
$L'$ and $N_q(X^{n-1})=\Theta_{n-1}^{\sqrt{q}+1,q}$, i.e.
$X^{n-1}\subseteq H'=\mathbb{P}^n$ is as in $(*)$ for any hyperplane $H'$ as above. 

For $i\neq 1$, consider now a hyperplane $H_i$ such that $L_{11}\subset
H_i$. Hence by Remark \ref{Remark3} we get $X^n\cap H_i=\cup_{j=1}^{\sqrt{q}+1}L_{ij}$,
where $L_{ij}=\mathbb{P}^{n-1}$ and $L_{i1}=L_{11}$. Define
$\Lambda_i:=\cap_{j=1}^{\sqrt{q}+1} L_{ij}$ and note that $\Lambda_i=\mathbb{P}^{n-2}$. 

If for any $H_i=\mathbb{P}^n$ as above such that $L_{11}\subset
H_i$ we have $\Lambda_i=\Lambda_1$ for $i=1,2,...,q+1$, then
$X^{n}=\cup_{i=1}^{q+1}\left(\cup_{j=1}^{\sqrt{q}+1} L_{ij}\right)$
with $\cap_{i,j}L_{ij}=\cap_{i=1}^{q+1}\Lambda_i=\Lambda_1=\mathbb{P}^{n-2}$. 
By taking a plane $P=\mathbb{P}^2$ such that $P\cap\Lambda_1=\emptyset$, we
see that $X^n$ is a cone with vertex $\Lambda_1=\mathbb{P}^{n-2}$
over a curve $\Gamma\subset P$ which is the
intersection set of all the $L_{ij}$'s with $P$, i.e. $\Gamma:=P\cap \left( \cup_{i,j}L_{ij}\right)$. 
Therefore, we get 
$$(d-1)q^n+dq^{n-1}+N_q(\mathbb{P}^{n-2})=N_q(X)=q^{n-1}N_q(\Gamma)+N_q(\mathbb{P}^{n-2}),$$
i.e. $N_q(\Gamma)=(d-1)q+d$, where $d=\sqrt{q}+1$.
Moreover, observe that $\Gamma\subset P=\mathbb{P}^2$ is a plane $(\sqrt{q}+1)$-arc. 
Thus from {\em \cite{B}} we know that $q\equiv 0\
(\mathrm{mod}\ d)$, i.e. $q=h(\sqrt{q}+1)$ for some
$h\in\mathbb{Z}_{\geq 1}$. Hence we get $h=p^t$ for some
$t\leq s$, where $q=p^s$, but this gives the numerical
contradiction $p^{s-t}=p^{\frac{s}{2}}+1$, since $p\geq 2$ and $s\geq 2$ even. 

This shows that there exists at least a hyperplane $\mathbb{P}^{n}$, say $H_2$, such that
$L_{11}\subset H_2$, $X^n\cap H_2=\cup_{j=1}^{\sqrt{q}+1}L_{2j}$, $L_{2j}=\mathbb{P}^{n-1}$,
$L_{21}=L_{11}$ and $\Lambda_2\neq\Lambda_1$. Set 
$W:=L_{11}\cap L_{12}\cap L'=\Lambda_1\cap L'=\mathbb{P}^{n-3}$ and define $W_{ij}:=V\cap L_{ij}$.

If for every $i,j$, we have $W_{ij}=W$, then $X^n=\cup_{i=1}^{q+1}\left( \cup_{j=1}^{\sqrt{q}+1} L_{ij}\right)$
and 
$$\mathbb{P}^{n-3}=W=\cap_{i,j}\left( W\cap L_{ij}\right) \subseteq\cap_{i,j}L_{ij}\subseteq (L_{11}\cap L_{12})\cap (L_{21}\cap L_{22})=\mathbb{P}^{n-3},$$ i.e. $W=\cap_{i,j}L_{ij}=\mathbb{P}^{n-3}$. By taking $U=\mathbb{P}^3$ such that
$W\cap U=\emptyset$, write $S':=U\cap\left( \cap_{i,j}L_{ij}\right)$. This shows that $X^n$ is a cone 
with vertex $W=\mathbb{P}^{n-3}$ over a surface $S'\subseteq U=\mathbb{P}^3$
without linear $\mathbb{F}_q$-components and such that $U\cap L_{ij}=\mathbb{P}^1\subset S'$ for every $i,j$.
Thus from Remark \ref{Remark4}(c) it follows that $N_q(S')=\Theta_{2}^{\sqrt{q}+1,q}:=(d-1)q^2+dq+1$,
i.e. $S'$ is a nonsingular Hermitian surface in $\mathbb{P}^3$ by \cite{HK5}. 

\smallskip

Therefore, without loss of generality,
we can assume that $W\neq W_{22}=W\cap L_{22}=\mathbb{P}^{n-4}$. So, after a change
of coordinates, we can write
$$L_{11}=L_{21}=V(X_1,X_2), \ L_{12}=V(X_1,X_2-X_n), $$
$$L_{22}=V(X_2,X_1-X_{n+1})\ \ \mathrm{and} \ \ L'=V(X_0,X_1).$$
Furthermore, since $L'\subseteq V(X_0)$, observe that $X^{n}\cap V(X_0)=X^{n-1}$ is as in $(*)$
by the inductive argument. So $X^n:=\{F=0\}$ is defined by the following homogeneous polynomial
$$F:=X_0f_0+(X_1,...,X_{n+1})\ A\ {}^t \! (X_1^{\sqrt{q}},...,X_{n+1}^{\sqrt{q}}),$$
where $f_0\in\mathbb{F}_q[X_0,...,X_{n+1}]$ and $A$ is an
$(n+1)\times (n+1)$ matrix such that ${}^t \! A=A^{(\sqrt{q})}$. Since
$H_1:=\langle L_{11}\cup L_{12} \rangle=\left\{X_1=0\right\}$,
$$X^n\cap H_1 =\left(\cup_{j=1}^{\sqrt{q}+1}L_{1j}\right)\ \mathrm{and}\ \ \Lambda_1=\{X_1=X_2=X_n=0 \}
= \cap_{j=1}^{\sqrt{q}+1}L_{1j},$$
we deduce that
$$X^n\cap\left\{X_1=0\right\}=\Pi_{i=1}^{\sqrt{q}+1}\left(a_2^iX_2+a_n^iX_n\right).$$ This shows that
$$F=X_0X_1g_0+(X_1,...,X_{n+1})\ A\ {}^t \! (X_1^{\sqrt{q}},...,X_{n+1}^{\sqrt{q}}),$$
where $g_0\in\mathbb{F}_q[X_0,...,X_{n+1}]$. Similarly, since
$H_2:=\langle L_{11}\cup L_{22} \rangle=\left\{X_2=0\right\}$,
$$X^n\cap H_2 =\left(\cup_{j=1}^{\sqrt{q}+1}L_{2j}\right)\ \mathrm{and}\
\Lambda_2=\{X_1=X_2=X_{n+1}=0 \}= \cap_{j=1}^{\sqrt{q}+1}L_{2j},$$ it follows that
$$X^n\cap\left\{X_2=0\right\}=\Pi_{i=1}^{\sqrt{q}+1}\left(a_1^iX_1+a_{n+1}^iX_{n+1}
\right).$$ This gives
$$F=X_0X_1X_2h_0+(X_1,...,X_{n+1})\ A\ {}^t \! (X_1^{\sqrt{q}},...,X_{n+1}^{\sqrt{q}}),$$
where $h_0\in\mathbb{F}_q[X_0,...,X_{n+1}]$. Furthermore, by an inductive
argument we know that $X^n\cap\left\{X_0-\lambda X_1=0\right\}$
is also a Hermitian variety for any $\lambda\in\mathbb{F}_q$. Thus
the polynomial
$$\bar{F}:=\lambda X_1^2X_2h_0(\lambda X_1,X_1,...,X_{n+1})$$ is Hermitian for any $\lambda\in\mathbb{F}_q^*$.
Since $\bar{F}$ does not contain the monomial $X_1X_2^{\sqrt{q}}$,
we deduce that $\lambda X_1^2X_2h_0(\lambda X_1,X_1,...,X_{n+1})$
is identically zero. In particular, $h_0(\lambda X_1,X_1,...,X_{n+1})$
is identically zero for any $\lambda\in\mathbb{F}_q^*$. Hence $(X_0-\lambda X_1)$ divides $h_0$ for
any $\lambda\in\mathbb{F}_q^*$. Since $\deg h_0=\sqrt{q}-2<q-1$,
we deduce that $h_0$ is the zero polynomial, i.e.
$$F=(X_1,...,X_{n+1})\ A\ {}^t \! (X_1^{\sqrt{q}},...,X_{n+1}^{\sqrt{q}}).$$
This shows that $X^n=\{F=0\}$ is a cone over a Hermitian variety $X^{n-1}$ for $n\geq
3$ and by an inductive argument we can conclude that
$X^n\subset\mathbb{P}^{n+1}$ is a cone over a nonsingular
Hermitian surface $X^2\subset\mathbb{P}^3$.
\hfill{$\square$}

\begin{prop}[$d=q+1$]\label{$d=q+1$}
Let $X^n\subset\mathbb{P}^{n+1}$ be a hypersurface of degree $q+1$
without $\mathbb{F}_q$-linear components and such that
$N_q(X^n)=\Theta_n^{q+1,q}$. If $n\geq 2$ then $X^n$ is a space filling hypersurface
$$(X_0,...,X_{n+1})\ A\ {}^t \! (X_0^q,...,X_{n+1}^q)=0,$$ where $A=\left( a_{ij}\right)_{i,j=1,...,n+2}$ is an
$(n+2)\times (n+2)$ matrix such that ${}^t \! A=-A$ and $a_{kk}=0$ for every $k=1,...,n+2$. Moreover, 
we have the following properties:
\begin{enumerate}
\item $X^n$ is singular if and only if $\det A = 0$;
\item if $n$ is odd, then $X^n$ is singular.
\end{enumerate}
\end{prop}

\noindent\textit{Proof}. Let $F=F(X_0,...,X_{n+1})$ be a
homogeneous polynomial of degree $q+1$ identically zero on
$\mathbb{F}_q^{n+2}$. From the Polynomial Evaluation Theorem {\em
\cite[Theorem 8]{C}} we know that the ideal of the algebraic set
$\mathbb{P}^{n+1}$ on $\mathbb{F}_q$ is generated by
$\left\{X_iX_j^q-X_i^qX_j\ |\ 0\leq i<j\leq n+1 \right\}.$ Thus
$F=\sum_{i,j=0}^{n+1}a_{ij}(X_iX_j^q-X_i^qX_j)$, i.e. $X^n$ is
defined by
$$(X_0,...,X_{n+1})\ A\ {}^t \! (X_0^q,...,X_{n+1}^q)=0,$$
where $A=\left( a_{ij}\right)_{i,j=1,...,n+2}$ is an
$(n+2)\times (n+2)$ matrix such that ${}^t \! A=-A$ and $a_{kk}=0$ for every $k=1,...,n+2$.
Write $F=\vec{x}\ A \ {}^t \! \vec{x}^{(q)}$, where $\vec{x}:=(X_0,...,X_{n+1})$ 
and $\vec{x}^{(q)}:=(X_0^q,...,X_{n+1}^q)$. Then $X^n$ is defined by $\vec{x}\ A \ {}^t \! \vec{x}^{(q)}=0 $ and we have
$$
\left(
\begin{array}{c}
\frac{\partial F}{\partial X_0} \\[3pt]
\vdots \\[3pt]
\frac{\partial F}{\partial X_{n+1}}
\end{array}\right)=A\cdot 
\left(
\begin{array}{c}
X_0^q \\[3pt]
\vdots \\[3pt]
X_{n+1}^q
\end{array}\right) =: A\cdot {}^t \! \vec{x}^{(q)} .
$$
Thus we see that $X^n$ is singular $\iff \det A=0$, and this proves $(1)$. 

Assume now that $n$ is odd. To prove $(2)$, by $(1)$ it is sufficient to show that $\det A=0$.
If $q$ is odd, then $\det A=\det {}^t \! A = \det (-A) = (-1)^{n+2}\det A=-\det A,$ i.e. $\det A=0$.  
Suppose now that $q$ is even and put $m:=n+2$. Then ${}^t \! A=A$ and from
$$sgn (\sigma)\left( \Pi_{i=1}^{m}a_{i,\sigma (i)}\right) =sgn (\sigma^{-1})\left( \Pi_{i=1}^{m}a_{\sigma (i),i}\right) =sgn (\sigma^{-1})\left( \Pi_{j=1}^{m}a_{j,\sigma^{-1} (j)}\right) $$
for any $\sigma\in S_m$, where $S_m$ is the symmetric group, we deduce that 
$$\det A = \sum_{\sigma\in S_m} sgn (\sigma)\left( \Pi_{i=1}^{m}a_{i,\sigma (i)}\right) =
\sum_{\sigma\in S_m :\ \sigma=\sigma^{-1}} sgn (\sigma)\left( \Pi_{i=1}^{m}a_{i,\sigma (i)}\right).$$
Note that any $\sigma\in S_m$ such that $\sigma=\sigma^{-1}$ is a product of distinct transpositions and since $m$ is odd, we deduce that
there exists $h\in\{1,...,m\}$ such that $\sigma (h)=h$ for every $\sigma\in S_m$ such that $\sigma=\sigma^{-1}$.
Since $a_{h,h}=0$ for every $h=1,...,n+2$, we conclude that $\det A=0$ and this proves $(2)$.
\hfill{$\square$}

\bigskip
\bigskip

Finally, let us observe that Propositions \ref{Prop3}, \ref{$d=q+1$} and
Theorem \ref{main theorem} give Corollary \ref{Cor1} of the
Introduction.

\bigskip
\bigskip

\bibliographystyle{elsarticle-num}
\bibliography{<your-bib-database>}

\end{document}